\documentclass[12pt]{amsart}

\RequirePackage{amsthm,amsmath,amssymb,natbib,verbatim}

\usepackage{geometry}
\theoremstyle{plain}

\newtheorem{theorem}{Theorem}[section]
\newtheorem{proposition}{Proposition}[section]

\newtheorem{assumption}{Assumption}[section]

\theoremstyle{definition}
\newtheorem{example}{Example}

\newtheorem{definition}{Definition}[section]

\newcommand{\nnr}{[0,\infty)}
\newcommand{\parmspace}{\boldsymbol{\vartheta}}

\newcommand{\R}{\mathbb{R}}
\newcommand{\X}{\mathcal{X}}
\newcommand{\kernels}{\mathcal{G}_{\nu}}

\newcommand{\statespace}{\mathcal{W}}
\newcommand{\A}{\mathcal{A}}
\newcommand{\B}{\mathcal{B}}
\newcommand{\C}{\mathcal{C}}
\newcommand{\D}{\mathcal{D}}

\newcommand{\Fnu}{\mathcal{F}_{\nu}}
\newcommand{\T}{\mathcal{T}}
\newcommand{\given}{ \vert }

\newcommand{\phee}{\varphi}
\newcommand{\Phee}{\varPhi}

\DeclareMathOperator{\E}{E}
\DeclareMathOperator{\risk}{\mathsf{R}}
\DeclareMathOperator{\ird}{IRD}

\providecommand{\norm}[1]{\lVert#1\rVert}

\begin{document}
\title[Evaluating Default Priors]{Evaluating Default Priors with a
  Generalization of Eaton's Markov Chain} 

\author[B.~P.~Shea and
G.~L.~Jones]{Brian P.~Shea and Galin L.~Jones \\ School of Statistics
  \\ University of Minnesota} \date{\today} 

\thanks{Brian Shea (bshea@stat.umn.edu) is a graduate student and
  Galin Jones (galin@umn.edu) is an Associate Professor.  Research
  supported by the National Science Foundation and the National
  Institutes of Health.  The authors thank Morris Eaton for many
  helpful conversations.}


\begin{abstract} 
  We consider evaluating improper priors in a formal Bayes setting
  according to the consequences of their use.  Let $\varPhi$ be a
  class of functions on the parameter space and consider estimating
  elements of $\varPhi$ under quadratic loss.  If the formal Bayes
  estimator of every function in $\varPhi$ is admissible, then the
  prior is strongly admissible with respect to $\varPhi$.  Eaton's
  method for establishing strong admissibility is based on studying
  the stability properties of a particular Markov chain associated
  with the inferential setting.  In previous work, this was handled
  differently depending upon whether $\varphi \in \varPhi$ was bounded
  or unbounded.  We introduce and study a new Markov chain which
  allows us to unify and generalize existing approaches while
  simultaneously broadening the scope of their potential
  applicability. To illustrate the method, we establish strong
  admissibility conditions when the model is a $p$-dimensional
  multivariate normal distribution with unknown mean vector $\theta$
  and the prior is of the form $\nu(\|\theta\|^{2})d\theta$.
\end{abstract}

\maketitle

\section{Introduction} 
\label{sec:intro}
Suppose we are in a parametric setting, and we are considering use of
an improper prior measure that yields a proper posterior distribution.
Such priors arise in the absence of honest prior belief about
parameter values and are typically derived from structural arguments
based on the likelihood or the parameter space \cite{KW:1996}.  Thus,
an improper prior, rather than being a statement of beliefs specific
to a situation, is a default. Such priors, proposed from likelihood or
invariance arguments, require evaluation, just as estimators proposed
from likelihood or invariance arguments require evaluation, and an
attractive avenue is to evaluate the prior according to the
consequences of its use.  That is, we can evaluate the prior by
examining properties of the resulting posterior inferences.  The
criterion we use to judge posterior inferences is known as
\textit{strong admissibility}.  This concept was introduced by Eaton
\cite{Eaton:1992} and has given rise to a substantial theory
\cite{Eaton:1992, Eaton:1997,Eaton:2001,Eaton:2004, FrenchAnnals:2007,
  Eaton:2008,Hobert:1999,hobe:schw:2002,hobe:tan:liu:2007,lai:1996}.
Our goal in the rest of this section is to convey the basic idea
behind strong admissibility and the way it is studied here.  We also
summarize our main results without delving too far into the details,
which are dealt with carefully later.

Suppose the sample space $\X$ is a Polish space with Borel
$\sigma$-algebra $\B$ and the parameter space $\parmspace$ is a Polish
space with Borel $\sigma$-algebra $\C$.  Let $\{P( \cdot | \theta),\;
\theta \in \parmspace\}$ be a family of sampling distributions where
we assume that for each $B \in \B$, $P(B|\cdot)$ is $\C$-measurable
and for each $\theta \in \parmspace$, $P(\cdot | \theta)$ is a
probability measure on $(\X, \B)$.  Let $\nu$ be a $\sigma$-finite
measure on the parameter space $\parmspace$ with
$\nu(\parmspace)=\infty$.  Throughout the marginal on $\X$
\begin{equation}
\label{eq:marginal}
M(dx) := \int_{\parmspace} P(dx | \theta) \, \nu(d\theta) \; .
\end{equation}
is assumed to be $\sigma$-finite.  In this case, the disintegration
\begin{equation}
\label{gen.bayes.eq} 
Q(d\theta \given x) \, M(dx) = P(dx \given \theta) \, \nu(d\theta)
\end{equation}
generalizes Bayes theorem and implicitly defines formal posterior
distributions $Q( \cdot \given x)$ on the parameter space.  Note that
for each $x \in \X$, $Q(\cdot | x)$ is a probability measure on
$(\parmspace, \C)$ and for each $C \in \C$, $Q(C |\cdot)$ is
$\B$-measurable. Taraldsen and Lindqvist \cite{tara:lind:2010} provide
a recent, accessible introduction to the existence of formal posterior
distributions while one can consult Eaton \cite{Eaton:1982,Eaton:1992}
and Johnson \cite{john:1991} for more details and references.

Suppose $\varphi : \parmspace \to \mathbb{R}^{p}$ for $p \ge 1$, and
consider estimating $\varphi(\theta)$.  The formal Bayes estimator of
$\varphi$ under squared error loss is the posterior mean
\begin{equation}
\label{eq:gen.bayes}
\hat{\varphi}(x) := \int_{\parmspace} \varphi(\theta) \, Q(d\theta | x)
\, .
\end{equation}
Let $\|\cdot\|$ denote the usual Euclidean norm. If $\delta(x)$ is any
estimator of $\phee(\theta)$, the risk function of $\delta$ is
\begin{equation}
\label{eq:risk}
\risk(\delta ; \theta) := \int_{\X} \norm{ \phee(\theta) - \delta(x)
}^2 P(dx \given \theta) \, .
\end{equation}
The estimator $\delta$ is \textit{almost-$\nu$ admissible} if for any
other estimator  $\delta_{1}$ such that  $\risk(\theta;\delta_{1}) \leq
\risk(\theta;\delta)$ for all $\theta \in \parmspace$, then the set
$\{ \theta \, : \, \risk(\theta;\delta_{1}) <
\risk(\theta;\delta) \}$ has $\nu$-measure 0.  

Since we will use admissibility to judge the prior, our interests are
more ambitious than establishing admissibility of a single
estimator. Let $\varPhi$ be a class of functions defined on the
parameter space.  If the formal Bayes estimator of every $\phee \in
\varPhi$ is almost-$\nu$ admissible, then we say the prior
(equivalently the posterior) is \textit{strongly admissible with
  respect to} $\varPhi$. A prior is strongly admissible if it is
robust against risk dominance within the class $\varPhi$. Since a
default prior will undergo repeated use, it is important for the range
of appropriate uses to be clearly defined and desirable that the range
be as large as possible.  We can then endorse the improper prior
insofar as it avoids unreasonable actions in a variety of such
problems.

Previous work on strong admissibility focused on the case where
$\varPhi$ consisted of a single unbounded function
\cite{Berger:1996,Berger:2005,Eaton:2001}, or all bounded functions
\cite{Eaton:1992, Eaton:1997,Eaton:2004, FrenchAnnals:2007,
  Eaton:2008,Hobert:1999,hobe:schw:2002,lai:1996}.  In either case, strong
admissibility was established either by verifying sufficient
conditions for the admissibility of an estimator established by
\citet{Brown:1971} or via Markov chain arguments using an approach
developed by \citet{Eaton:1992,Eaton:2001}.  We study the latter
method.

Eaton's method for establishing almost-$\nu$ admissibility of formal
Bayes estimators is based on the recurrence properties of a Markov
chain associated with the inferential setting; the relevant notion of
recurrence is defined in the next section.  However, different Markov
chains were required depending upon whether $\varphi$ was bounded
\cite{Eaton:1992} or unbounded \cite{Eaton:2001}. We introduce and
study a new Markov chain which allows us to unify and generalize these
existing approaches while simultaneously broadening the scope of their
potential applicability. The expected posterior
\begin{equation}
\label{Eaton.kernel}
R(d\theta \given \eta) = \int_{\X} Q(d\theta \given x) \, P(dx \given \eta) \; . 
\end{equation}
is a Markov kernel on $\parmspace$.  We study transformations of $R$,
which are now described.  Let $f : \parmspace \times \parmspace \to
\mathbb{R}^{+}$ satisfy $f(\theta, \eta) = f(\eta, \theta)$ and set
\[
\mathcal{T}(\eta) = \int_{\parmspace} f(\theta, \eta) R(d\theta |
\eta) \; .
\]
If $0 < \mathcal{T}(\eta) < \infty$ for all $\eta
\in \parmspace$, then 
\[
T(d\theta | \eta) = \frac{f(\theta, \eta)}{\mathcal{T}(\eta)}
R(d\theta | \eta) 
\]
is a Markov kernel on $\parmspace$. Recurrence of the Markov chain
associated with $T$ implies the almost-$\nu$ admissibility of formal
Bayes estimators with respect to a large class of functions.  Define
\begin{equation}
\label{eq:Phee f}
\varPhi_{f} := \left\{ \varphi \, : \,  \|\varphi(\theta) -
  \varphi(\eta) \|^{2} \le M_{\varphi} f(\theta, \eta) \; \text{ some
  } 0 < M_{\varphi} < \infty \right\}\; . 
\end{equation}
We prove that if the Markov chain defined by $T$ is recurrent, then
the formal Bayes estimator of every function in $\varPhi_{f}$ is
almost-$\nu$ admissible, and we say the prior $\nu$ is
\textit{strongly admissible with respect to $\varPhi_{f}$}.  The
following example illustrates this technique.

\begin{example}
\label{ex:normal-lebesgue}
Let $X$ be a $p$-dimensional normal random variable with identity
covariance $I_{p}$ and unknown location $\theta$.  Let $p$-dimensional
Lebesgue measure $\nu(d\theta) = d\theta$ be our improper prior. The
proper posterior for $\theta$ is a normal with mean $x$---the observed
value of the random variable $X$---and covariance $I_{p}$. The kernel 
\[
R(d\theta \given \eta)  = \int_{\X} Q(d\theta \given x) P(dx \given
\eta) =  (4\pi)^{-p/2} \exp \left( - \frac{1}{4} \lVert \theta - \eta
  \rVert^{2} \right) d\theta  
\]
describes a $N(\eta,2I_{p})$ random variable.  Let $d$ be an arbitrary
positive constant and note that $f(\theta, \eta ) = \lVert \theta -
\eta \rVert^{2} + d$ is symmetric and uniformly bounded away from
0. Further
\[ 
\mathcal{T}(\eta) = \int f(\theta, \eta)  R(d\theta \given \eta) = \int\left( \lVert \theta - \eta \rVert^{2} + d \right) R(d\theta \given \eta) = 2p +d. 
\]
Thus,
\[ 
T(d\theta \given \eta) = \frac{\lVert \theta - \eta \rVert^{2} + d}{2p + d} R(d\theta \given \eta) 
\] 
and the chain with kernel $T$ is a random walk on $\R^{p}$. Note that 
\[
\int \lVert \theta - \eta \rVert^{p} \, T(d\theta \given \eta) <
\infty \] since a normal distribution  has moments of all orders and
hence for $p=1$ or $p=2$, the chain is recurrent \cite{Chung:1951,
  Revuz:1984}.  We conclude that  for $p=1$ or $p=2$ Lebesgue measure
is strongly admissible.  That is, the formal Bayes estimators of
all functions $\varphi$ satisfying
\[
\|\varphi(\theta) - \varphi(\eta)\| \le M_{\varphi}\left( \|\theta -
  \eta\|^{2} + d \right) \; \text{ for some } \; 0 < M_{\varphi} < \infty \; 
\]
are almost-$\nu$ admissible.  Since $d$ is arbitrary this includes
all bounded functions as well as many unbounded functions.

Lebesgue measure is not strongly admissible with respect to
$\varPhi_{f}$ when $p \ge 3$ since the James-Stein estimator dominates
the formal Bayes estimator of $\phee(\theta)=\theta$.
In Section~\ref{sec:normal} we consider the normal means problem with
an alternative prior when $p \ge 3$.  We will return to this example
below.
\end{example}

Recurrence of the Markov chain described by $T$ implies more than we
have so far claimed.  Suppose the function $u \colon \parmspace \to
\mathbb{R}^{+}$ is bounded away from zero and infinity so that $1/c <
u(\theta) < c$ for some constant $c > 1$ and every $\theta
\in \parmspace$. The measure defined by
\begin{equation}
\label{eq:bounded perturbation}
\nu_{u}(d\theta) :=  u(\theta) \,\nu(d\theta) 
\end{equation}
is a \textit{bounded perturbation} of $\nu$.  Let $\Fnu$ be the family
of all bounded perturbations of $\nu$.  Observe that $\nu \in \Fnu$
and that the other elements of $\Fnu$ are measures with tail behavior
similar to $\nu$.  Eaton \cite{Eaton:1992} showed that recurrence of
the chain with kernel $R$ implied the formal Bayes estimators of any
bounded function is almost-$\nu$ admissible for every prior in $\Fnu$.
We extend this result and show that recurrence under $T$ is sufficient
for the strong admissibility with respect to $\varPhi_{f}$ of every
element of $\Fnu$.  In this case, we say the family $\Fnu$ is
\textit{strongly admissible with respect to $\varPhi_{f}$}.

\begin{example} \label{ex:normal-lebesgue2}
Recall the setting of Example~\ref{ex:normal-lebesgue}. Let $\Fnu$ be
the bounded perturbations of $p$-dimensional Lebesgue measure
$\nu(d\theta)=d\theta$.  That is, elements of $\Fnu$ are measures of
the form $\nu_{u}(d\theta) = u(\theta) d\theta$ where $1/c < u(\theta)
< c$ for some constant $c > 1$.  Recall that $\varPhi_{f}$ is the
class of all bounded functions and all functions satisfying 
\[
\|\varphi(\theta) - \varphi(\eta)\| \le M_{\varphi}\left( \|\theta -
  \eta\|^{2} + d \right) \; \text{ for some } 0 < M_{\varphi} < \infty
\; . 
\]  
The recurrence of the Markov chain governed by $T$ when $p=1$ or $p=2$
implies the formal Bayes estimators of every $\varphi \in \varPhi_{f}$
are almost-$\nu_{u}$ admissible for every prior $\nu_{u} \in \Fnu$.
That is, $\Fnu$ is strongly admissible with respect to $\varPhi_{f}$
 \end{example}

 Our main results generalize existing work \cite{Eaton:1992,
   Eaton:2001} and in fact unify the analysis for bounded and
 unbounded functions.  We show that the Markov kernel $R$ can be
 transformed to define many Markov chains, any one of which might be
 used to demonstrate strong admissibility, thus greatly broadening the
 scope of potential applications.  We are never concerned narrowly
 with a single admissibility problem but broadly with species of
 problems.  Moreover, solving a single representative problem, which
 representative we are free to elect, solves any problem within a
 bounded rate of change---whether of the function to be estimated or
 the prior used to estimate it.

 The remainder is organized as follows.
 Section~\ref{sec:preliminaries} gives some background on recurrence
 for general state space Markov chains.
 Section~\ref{sec:radon-nikodym} presents the main results and
 Section~\ref{sec:normal} illustrates the main results by considering
 the multivariate normal means problem.  Finally, many technical
 details are deferred to the appendices.
 
\section{Recurrence of Markov chains} \label{sec:preliminaries} 

The goal of this section is to introduce a general notion of
recurrence for Markov chains. Let $\statespace$ be a Polish space and
denote the Borel $\sigma$-algebra by $\D$.  Let $K : \D \times
\statespace \to [0,1]$.  Then $K$ is a \textit{Markov transition
  kernel} on the measurable space $(\statespace,\D)$ if $K(D \given
\cdot)$ is a nonnegative measurable function for every $D \in \D$ and
$K(\cdot \given w)$ is a probability measure for every $w \in
\statespace$.

The kernel $K$ determines a time-homogeneous Markov chain $W=\{ W_{0},
W_{1}, W_{2}, \ldots\}$ on the product space $\statespace^{\infty}$
which is equipped with the product $\sigma$-algebra $\D^{\infty}$.
Note that conditional on $W_{n}=w$, the law of $W_{n+1}$ is $K(\cdot
\given w)$.  Given $W_{0} = w$, let $\Pr(\cdot \given w)$ be the law
of $W$ on $\statespace^{\infty}$.

Suppose $D \in \D$.  The random variable
\[ 
\tau_{D} = \begin{cases} \infty \text{ if $W_{n} \notin D$ for all
    $n\geq 1$, or} \\ \text{the smallest $n\geq 1$ such that $W_{n} \in
    D$ otherwise}  \end{cases}  
\]
is a stopping time for $D$, and $E_{D} = \{ \tau_{D} < \infty \}$ is
the set of paths that encounter $D$ after initialization. Let $\xi$ be
a non-trivial, $\sigma$-finite measure on $(\statespace,\D)$ and
recall that a set is \textit{$\xi$-proper} if its measure under $\xi$
is positive and finite.

\begin{definition}
  A $\xi$-proper set $D$ is \textit{locally $\xi$-recurrent} if
  $\Pr(E_{D} \given w) = 1$ for all but a $\xi$-null set of
  initial values in $D$.  Call the Markov chain $W$ \textit{locally
  $\xi$-recurrent} if every $\xi$-proper set is locally
  $\xi$-recurrent.
\end{definition}

This notion of recurrence is more general than that typically
encountered in general state space Markov chain theory
\cite{meyn:twee:1993}, but is appropriate since the chains we will
consider in the next section may not be irreducible
\cite{FrenchAnnals:2007, hobe:tan:liu:2007}.  We consider a method for
establishing local recurrence in Section~\ref{sec:reducing dimension}.

The connection between the above general Markov chain theory and the
notion of strong admissibility relies heavily on the special structure
of symmetric Markov chains.  Let $\xi$ be a non-trivial,
$\sigma$-finite measure.  Then the kernel $K$ is $\xi$-symmetric if
the measure
\[
\omega(A,B) : = \int_{A} K(B|w) \xi(dw)
\] 
satisfies $\omega(A, B) = \omega(B, A)$ for all $A, B \in \D$.
Throughout the remainder we restrict attention to symmetric Markov
kernels.  Eaton \cite{Eaton:1997,Eaton:2004} provides some background
on the theory of symmetric Markov chains underlying strong
admissibility.

\section{Strong Admissibility via Markov
  chains} \label{sec:radon-nikodym} 

In a ground-breaking paper Eaton \cite{Eaton:1992} connected the local
recurrence of a Markov chain and the almost-$\nu$ admissibility of the
formal Bayes estimator for bounded $\phee$.  Later, Eaton
\cite{Eaton:2001} showed that the local recurrence of a different
Markov chain was required to establish almost-$\nu$ admissibility when
$\phee$ is unbounded. Since our work builds on them, these
foundational definitions and results are stated carefully here.  The
expected posterior distribution at the parameter value $\eta$,
$R(d\theta | \eta)$ defined in \eqref{Eaton.kernel}, is a
$\nu$-symmetric Markov kernel on $\parmspace$ since it satisfies the
detailed balance condition
\begin{equation}
R(d\theta \given \eta) \, \nu(d\eta) = R(d\eta \given \theta) \,
\nu(d\theta) \; . \label{detailed.balance} 
\end{equation}
Call $R$ an \textit{Eaton kernel}--analogous kernels exist for the
sample space \citep{Hobert:1999} and the product of the sample and
parameter spaces \citep{FrenchAnnals:2007}.  Eaton \cite{Eaton:1992}
established the following basic result. 
\begin{theorem}
\label{thm:eaton bdd}
If the Markov chain with kernel $R$ is locally-$\nu$ recurrent, then
the formal Bayes estimator of every bounded function is almost-$\nu$
admissible.
\end{theorem}

Eaton \cite{Eaton:1992} also showed that the result holds for bounded
perturbations of $\nu$.  Theorem~\ref{thm:eaton bdd}  has found
 substantial application \cite{Eaton:1992, Eaton:2004,
  FrenchAnnals:2007, Eaton:2008,Hobert:1999,hobe:schw:2002,lai:1996}. 

Now suppose we want to use Markov chains to study formal Bayes
estimators of unbounded functions on $\parmspace$.  We need a
basic assumption on the risk function $\risk$, defined at
\eqref{eq:risk}, to ensure existence of the integrated risk
difference which is studied in the appendices.
\begin{assumption}
\label{ass:risk}
Suppose there exist sets $K_{1} \subseteq K_{2}
\subseteq \cdots$ with $\cup K_{i} = \parmspace$.  Also assume that
for all $i$, $0 < \nu(K_{i}) < \infty$ and
\[
\int_{K_{i}} \risk(\hat{\varphi}; \theta) \nu(d\theta) < \infty \; .
\]
\end{assumption}
Let $\phee : \parmspace \to \mathbb{R}^{p}$ be measurable and
unbounded and set
\[
\mathcal{S}(\eta) = \int_{\parmspace}  \lVert \phee(\theta) -
\phee(\eta) \rVert^{2} \, R(d\theta \vert \eta)  
\] 
so that if $0 < \mathcal{S}(\eta) < \infty$ for all $\eta$, then
\begin{equation}
\label{Eaton.kernel.ubdd}
S(d\theta | \eta) = \frac{\lVert \phee(\theta) - \phee(\eta)
  \rVert^{2}}{\mathcal{S}(\eta)} \, R(d\theta \vert \eta)  
\end{equation}
is a Markov kernel defining a Markov chain on $\parmspace$.  Moreover $S$
satisfies detailed balance with respect to
$\mathcal{S}(\eta)\nu(d\eta)$: 
\[ 
S(d\theta|\eta) \mathcal{S}(\eta)\nu(d\eta) = S(d\eta | \theta)
\mathcal{S}(\theta) \nu(d\theta) \; .   
\] 
Eaton \cite{Eaton:2001} connected the local recurrence of $S$ with the
almost-$\nu$ admissibility of the formal Bayes estimator of
$\phee(\theta)$. 

\begin{theorem}\label{thm:eaton ubdd} 
  Suppose Assumption~\ref{ass:risk} holds and assume $0 <
  \mathcal{S}(\eta) < \infty$ for each $\eta \in \parmspace$ and there
  exist $\nu$-proper sets $K_{1} \subseteq K_{2} \subseteq \cdots$
  such that $\cup K_{i} = \parmspace$, and for each $i$
\[ 
\int_{K_{i}} \mathcal{S}(\eta) \nu(d\eta) < \infty \; .  
\] 
If the Markov chain with kernel $S$ is locally-$\nu$ recurrent, then
the formal Bayes estimator of $\phee$ is almost-$\nu$ admissible.  
\end{theorem}

In the next section, we unify and generalize these results. We show
that by analyzing an appropriate Markov chain we can recover the
conclusions of both theorems and, in fact, achieve something stronger.
Moreover, we broaden the class of Markov chains that can be studied to
obtain strong admissibility results.  

\subsection{A new Markov chain connection}

For any $\eta \in \parmspace$, let $\psi( \cdot \given \eta)$ be a
nontrivial $\sigma$-finite measure on $(\parmspace,\C)$ such that
$\psi( \cdot \given \eta)$ is absolutely continuous with respect to
$R( \cdot \given \eta)$.  For any element $C$ of the Borel sets $\C$,
let $\psi(C \given \cdot)$ be a nonnegative measurable function.  Let
$f( \cdot, \eta)$ be a Radon-Nikodym derivative of $\psi( \cdot \given
\eta)$ with respect to $R( \cdot \given \eta)$---  
\begin{equation}
 \psi(C \given \eta) = \int_{C} f(\theta, \eta) \, R(d\theta \given
 \eta) \quad \text{for all} \, (\eta,C) \in (\parmspace,\C)
 . \label{psi.def} 
\end{equation} 
Define
\begin{equation}
T(C \given \eta) = \frac{\psi(C \given \eta)}{\psi(\parmspace \given
  \eta)} \quad \text{for all} \, (\eta,C) \in
(\parmspace,\C). \label{T.def} 
\end{equation}
Let $\mathcal{T}(\eta)=\psi(\parmspace \given \eta)$ for all $\eta
\in \parmspace$.  We will make the following basic assumptions on
$\mathcal{T}$. 
\begin{assumption}
\label{ass:T}
For all $\eta \in \parmspace$ we have $0 < \mathcal{T}(\eta) < \infty$
and there exist $\nu$-proper sets $K_{1} \subseteq K_{2} \subseteq
\cdots$ such that $\cup K_{i} = \parmspace$, and for each $i$
\[ 
\int_{K_{i}} \mathcal{T}(\eta) \nu(d\eta) < \infty \; .  
\] 
\end{assumption}
The first part of the assumption ensures that the kernel $T$ is
well-defined. If $\mu(d\eta) = \mathcal{T}(\eta) \nu(d\eta)$, then the
second part of the assumption  implies that $\mu$ is $\sigma$-finite.
In Proposition~\ref{rnd.thm.1} we establish that 
\[
T(d\theta | \eta) \mu(d\eta) = f(\theta, \eta) R(d\theta | \eta) \nu(d\eta) \; .
\]
If $f(\theta, \eta) = f(\eta, \theta)$ for all $\eta$ and $\theta$ in
$\parmspace$, then  by using \eqref{detailed.balance} it is easy to
see that $T$ is symmetric with respect to $\mu$. 

Note that $\mathcal{T}(\eta)$ is almost-$\nu$ uniformly bounded away
from 0 if there exists $\epsilon > 0$ such that $\mathcal{T}(\eta) \ge
\epsilon$ except possibly on a set of $\nu$-measure 0.  We are now in
a position to state the main result.  The proof is given in
Appendix~\ref{app:main result}.  Recall the definition of $\Phi_{f}$
from \eqref{eq:Phee f}.

\begin{theorem} \label{thm:main result} Suppose
  Assumptions~\ref{ass:risk} and~\ref{ass:T} hold. Let
  $\mathcal{T}(\eta)$ be almost-$\nu$ uniformly bounded away from zero
  and suppose that for all $\eta$ and $\theta$ in $\parmspace$
\begin{equation}
f(\theta, \eta) = f(\eta, \theta) \; . \label{density.symmetry}
\end{equation}
If the Markov chain having kernel $T(d\theta | \eta)$ is locally
$\nu$-recurrent,  then $\nu$ is strongly admissible with respect to
$\Phee_{f}$. 
\end{theorem}

The following extension of Theorem~\ref{thm:main result} is also
proved in Appendix~\ref{app:main result}.

\begin{theorem}
\label{thm:family}
Assume the conditions of Theorem~\ref{thm:main result}.  Then every
bounded perturbation $\nu_{u}$ is strongly admissible, that is, the
family $\Fnu$ is strongly admissible with respect to $\Phee_{f}$. 
\end{theorem}

If we take $f(\theta, \eta) = d$ for some $d > 0$, then we completely
recover the results of Theorem~\ref{thm:eaton bdd} while if $f(\theta,
\eta) = \|\varphi(\theta) - \varphi(\eta)\|^{2}$, then we extend the
results of Theorem~\ref{thm:eaton ubdd}.  If, as in
Examples~\ref{ex:normal-lebesgue} and~\ref{ex:normal-lebesgue2}, we
set $f(\theta, \eta) = \|\varphi(\theta) - \varphi(\eta)\|^{2} + d$
with $d >0$, then we obtain results stronger than if we had
established local recurrence of the chains associated with the kernels
$R$ and $S$ and relied on Theorems~\ref{thm:eaton bdd} and
~\ref{thm:eaton ubdd}. Moreover, since the analyst has the freedom to
choose an appropriate $f$, this result extends the range of potential
applicability of Eaton's method.

\subsection{Reducing dimension}
\label{sec:reducing dimension}

The Markov kernel $T$ naturally takes the same dimension as the
parameter space.  This dimension may be quite large, making the
required analysis difficult.  In this section we prove that the
conclusions of Theorem~\ref{thm:family} (hence Theorem~\ref{thm:main
  result}) hold if we can establish the local recurrence of a
particular Markov chain which lives on $[0,\infty)$.

Denote the Borel subsets of $[0,\infty)$ by $\A$.  A measurable
mapping $t$ from $(\parmspace,\C)$ to $([0,\infty),\A)$ induces a
measure $\tilde{\nu}$ given by
\[
\tilde{\nu}(A) = \nu(t^{-1}(A)) \qquad A \in \A \; .
\]
We will need the next assumption throughout the remainder of this section.

\begin{assumption} \label{ass:partition}
There exists a partition $\{A_{i}\}$ of $[0,\infty)$ such that each
$A_{i}$ is measurable and each $C_{i} = t^{-1}(A_{i})$ is $\nu$-proper.
\end{assumption}

\citet{Eaton:2008} showed that under Assumption~\ref{ass:partition}
there exists a Markov transition function $\pi(d\theta \given \beta )$
on $\C \times [0,\infty)$ such that  
\begin{equation}
\nu(d\theta) = \pi(d\theta \given \beta) \tilde{\nu}(d \beta) \label{nu.tilde}
\end{equation}
which means for all measurable nonnegative functions $f_{1}$ on $\parmspace$
and $f_{2}$ on $[0,\infty)$ 
\[ 
\int_{\parmspace} f_{2}(t(\theta)) f_{1}(\theta) \, \nu(d\theta) =
\int_{0}^{\infty} f_{2}(\beta) \left( \int_{\parmspace} f_{1}(\theta) \,
  \pi(d\theta \given \beta) \right) \, \tilde{\nu}(d \beta) \; . 
\]
Define
\begin{equation} \label{p.tilde} 
\tilde{P}(dx \given \beta) = \int_{\parmspace} P(dx \given \theta)
\pi(d\theta \given \beta) \qquad \beta \in [0, \infty) \; . 
\end{equation}
The conditional probabilities $\tilde{P}( \cdot \given \beta)$ form a
parametric family indexed by $\beta $, and $\tilde{\nu}(d \beta)$ is a
$\sigma$-finite prior.  \citet{Eaton:2008} also showed that the
marginal measure on $\X$ is the same as at \eqref{eq:marginal}. That
is,
\[ 
\int_{0}^{\infty} \tilde{P}(dx \given \beta) \tilde{\nu}(d \beta) =  M(dx) , 
\] 
which is assumed $\sigma$-finite. Thus, there is a Markov kernel
$\tilde{Q}(d\beta \given x)$ satisfying 
\begin{equation}
\label{reduced.gen.bayes.eq}
\tilde{Q}(d\beta \given x) M(dx) = \tilde{P}(dx \given \beta)
\tilde{\nu}(d \beta) \, . 
\end{equation}
In fact,  a version of the posterior \cite{Eaton:2008} is 
\begin{equation}
\tilde{Q}(A \given x) := Q(t^{-1}(A) \given x) \quad \text{for all $A
  \in \A$} \; . \label{q.tilde} 
\end{equation}
The expected posterior
\begin{equation}
\tilde{R}(d\beta \given \alpha) = \int \tilde{Q}(d\beta \given x)
\tilde{P}(dx \given \alpha) \label{reduced.Eaton.kernel} 
\end{equation}
is a $\tilde{\nu}$-symmetric Eaton kernel.

If $\eta$ and $\theta$ are elements of $\parmspace$, let
$\alpha=t(\eta)$ and $\beta=t(\theta)$.  For any $\alpha \in \nnr$,
let $\tilde{\psi}( \cdot \given \alpha)$ be a nontrivial
$\sigma$-finite measure on $(\nnr,\A)$ such that $\tilde{\psi}( \cdot
\given \alpha)$ is absolutely continuous with respect to $\tilde{R}(
\cdot \given \alpha)$. Let $\tilde{f}( \cdot, \alpha)$ be a
Radon-Nikodym derivative of $\tilde{\psi}( \cdot \given \alpha)$ with
respect to $\tilde{R}( \cdot \given \alpha)$---ie, for all $(\alpha,A)
\in (\nnr,\A)$ 
\[ 
\tilde{\psi}(A \given \alpha) = \int_{A} \tilde{f}(\beta, \alpha) \, \tilde{R}(d\beta \given \alpha) \; .
\]
Define $\tilde{\T}(\alpha) = \tilde{\psi}([0, \infty) \given \alpha)$ and set
\[
\tilde{T}(A \given \alpha) = \frac{\tilde{\psi}(A \given \alpha)}{\tilde{\T}(\alpha)} \; \;  \qquad \; ~~~ (\alpha,A) \in (\nnr,\A)
\]
and let $\tilde{\mu}(d\alpha) = \tilde{\T}( \alpha) \,
\tilde{\nu}(d\alpha)$.  The following assumption ensures that
$\tilde{T}$ is a well-defined kernel and that $\tilde{\mu}$ is
$\sigma$-finite.
\begin{assumption}
\label{ass:tildeT}
For all $\alpha$ we have $0 < \tilde{\T}(\alpha) < \infty$ and there
exist $\tilde{\nu}$-proper sets $K_{1} \subseteq K_{2} \subseteq
\cdots$ such that $\cup K_{i} = [0, \infty)$, and for each $i$
\[ 
\int_{K_{i}} \tilde{\T}(\alpha) \tilde{\nu}(d\alpha) < \infty \; .  
\] 
\end{assumption}
In Appendix~\ref{app:reduced result} we show that $\tilde{T}$ is
$\tilde{\mu}$-symmetric. The following theorem shows that we can
analyze the recurrence properties of the chain defined by $\tilde{T}$
to achieve the conclusions of Theorem~\ref{thm:family} (hence
Theorem~\ref{thm:main result}).  The proof is given in
Appendix~\ref{app:reduced result} while use of the result is
illustrated in Section~\ref{sec:normal}.  Recall that $\Fnu$ is the
family of bounded perturbations of $\nu$.

\begin{theorem}\label{thm:reduced result}
Suppose Assumption~\ref{ass:tildeT} holds.  Let $\tilde{\T}(\cdot)$ be
almost-$\tilde{\nu}$ uniformly bounded away from 0 and suppose that
for all $\alpha$ and $\beta$ in $[0,\infty)$ 
\[
\tilde{f}(\beta, \alpha) = \tilde{f}(\alpha, \beta) \; .
\]
If the Markov chain with kernel $\tilde{T}$ is locally
$\tilde{\nu}$-recurrent, then $\Fnu$ is strongly admissible with
respect to $\Phi_{f}$. 
\end{theorem}

Since $\tilde{T}$ lives on $\nnr$ it would be convenient to have
conditions which guarantee the local recurrence of a Markov chain on
$\nnr$.  This is discussed in the following section.

\subsubsection{Recurrence of Markov chains on $\nnr$}
\label{sec:nnr}

To this point we have said little about establishing local recurrence.
The following theorem presents one method for doing so and is a
distillation of several existing results \cite{Eaton:2004,Eaton:2008}.
It applies generally to Markov chains on $\nnr$ and hence the notation
in this section is consistent with that of
Section~\ref{sec:preliminaries}. 

Let $\statespace = [0, \infty)$ and $\D$ be the Borel
$\sigma$-algebra.  Let $K : \D \times \statespace \to [0,1]$ be a
Markov kernel which defines a time-homogeneous Markov chain $W=\{
W_{0}, W_{1}, W_{2}, \ldots\}$ on $\statespace^{\infty}$.  
Define the $k$th moment of $K$ about its current state as
\begin{equation} 
\label{transition.moments} 
m_{k}(v) = \int_{0}^{\infty} (x-v)^{k} K(dx \given v) \, . 
\end{equation}

\begin{theorem} \label{thm:recurrence}
 Assume for each positive integer $n$ there exists $\delta(n) < 1$ such that 
\begin{equation}
\label{sup.T}
\sup_{v \in [0,n]} K( [0,n] \given v ) \leq \delta(n) \; . 
\end{equation}
Suppose
\begin{equation}
 \label{moment.cond.1}
\lim_{v \to \infty} \frac{\log v}{v} \frac{m_{3}(v)}{m_{2}(v)} = 0
\end{equation}
and there exists a function $\phi$ and an integer $n_0$ such that for
$v \in [n_0, \infty)$ 
\begin{equation}
\label{moment.cond.2}
m_{1}(v) \leq \frac{m_{2}(v)}{2v} [1+\phi(v)] \quad \text{and} \quad
\lim_{v \to \infty}\phi(v)\log v = 0 \, . 
\end{equation}
If $\xi$ is a non-trivial, $\sigma$-finite measure, $0 < \xi([0,n_0))
< \infty$ and $K$ is $\xi$-symmetric, then the Markov chain $W$ is
locally $\xi$-recurrent.
\end{theorem}

We will use Theorem~\ref{thm:recurrence} in conjunction with  Theorem~\ref{thm:reduced result}  in our main application in Section~\ref{sec:normal}.

\section{Admissible Priors for the Multivariate Normal
  Mean} \label{sec:normal} 

Let $X \sim N_{p}(\theta,I_{p})$ with $p \ge 1$. Consider the family
of $\sigma$-finite measures on $\parmspace = \R^{p}$ described by 
\begin{equation}
 \left( \frac{1}{a+\lVert \theta \rVert^{2}} \right)^{b}
 d\theta \label{priors} \quad \quad a \ge 0, \quad b \ge 0 \; . 
\end{equation}
For $a=0$, the prior is improper for all $b > 0$, but the induced
marginal distributions on $\X$ are only $\sigma$-finite for $b < p/2$.
For $a>0$, the family yields improper prior distributions when $b \leq
p/2$ and proper prior distributions when $b>p/2$.  In fact, if
$b=(a+p)/2$, the prior is the kernel of a multivariate $t$
distribution with $a$ degrees of freedom.  When $a>0$ and $b=0$, this
is $p$-dimensional Lebesgue measure which was considered in
Example~\ref{ex:normal-lebesgue}. 

Now suppose $p \ge 3$.  \citet{Berger:2005} established that the
formal Bayes estimator of $\theta$ is admissible when $a=1$ and $b=
(p-1)/2$ while \citet{Eaton:2008} use Theorem~\ref{thm:eaton bdd} to
prove that if $a >0$ and $b \in [p/2 -1 , \, p/2]$, then the formal
Bayes estimator of every bounded function is almost admissible. We use
Theorems~\ref{thm:reduced result} and~\ref{thm:recurrence} to add to
these results.  Let $d > 0$ be arbitrary and define
\[
f(\theta, \eta) = \|\theta - \eta\|^{2} + d \quad \quad \theta, \eta
\in \mathbb{R}^{p}  
\]
so that
\[
\varPhi_{f} = \{ \varphi \, : \, \|\varphi(\theta) - \varphi(\eta) \|
\le M_{\varphi} f(\theta, \eta) \quad \text{some} \; \; 0 <
M_{\varphi} < \infty \} \; . 
\]

\begin{theorem} \label{illustration}
For $p \ge 3$ let $X \sim N_{p}(\theta,I_{p})$ and with $a>0$ set 
\[ 
\nu(d\theta) = \left( \frac{1}{a+\lVert \theta \rVert^{2}} \right)^{p/2} d\theta \; .
\]
Then $\Fnu$, the family of bounded perturbations of $\nu$, is strongly
admissible with respect to $\Phee_{f}$.
\end{theorem}

\begin{proof}[Proof of Theorem \ref{illustration}]
Let $\gamma_{a}(z) = (a+z)^{-p/2}$, so that our family of priors \eqref{priors} can be expressed as $\gamma_{a}(\lVert \theta \rVert^{2}) \, d\theta$. The function $t(\theta) = \lVert \theta \rVert^{2}$ fulfills the requirements of Assumption \ref{ass:partition} as can be seen by letting the sets $A_{i} = [i-1,i)$ partition the nonnegative real numbers. Letting $\pi(d\theta \given \beta)$ denote the uniform distribution on the hypersphere of radius $\sqrt{\beta}$, it can be shown that $\nu(d\theta) = \pi(d\theta \given \beta) \tilde{\nu}(d\beta)$ with
\[
\tilde{\nu}(d\beta) = \frac{[\Gamma(1/2)]^{p}}{\Gamma(p/2)} \, \gamma_{a}(\beta) \, \beta^{\frac{p}{2}-1} \, d\beta \label{reduced.prior}
\]
on $\nnr$.  Let $\pi_{1}$ be the uniform distribution on the unit
hypersphere $\Xi$.  Then the reduced sampling distribution has density
\[ 
\tilde{p}(x \given \beta) = \int_{\Xi} (2\pi)^{-p/2} \exp\left( -\frac{1}{2} \lVert x - \xi \sqrt{\beta} \rVert \right) \, \pi_{1}(d\xi) 
\]
with respect to Lebesgue measure on $\X$. If 
\[
m(x) = \int_{0}^{\infty} \tilde{p}(x|\beta) \tilde{\nu}(d\beta),
\]
then the formal posterior has density 
\[ 
\tilde{q}( \beta \given x ) =  \frac{[\Gamma(1/2)]^{p}}{\Gamma(p/2)} \,\frac{ \tilde{p}(x \given \beta) \,\gamma_{a}(\beta) \, \beta^{\frac{p}{2}-1} }{m(x)} 
\]
with respect to Lebesgue measure on $\nnr$. The expected posterior
\begin{equation}
\tilde{R}(d\beta \given \alpha) = \int_{\X} \tilde{q}( \beta \given x ) \, \tilde{p}(x \given \alpha) \, dx \, d\beta
\end{equation}
is a $\tilde{\nu}$-symmetric Markov transition kernel on $[0,\infty)$.

Let $\tilde{f}(\beta, \alpha) = 2 (\beta + \alpha + c)$ with $c$ a positive constant.  It is clear that $\tilde{f}$ is symmetric in $\alpha$ and $\beta$ and bounded away from zero. For any nonnegative real number $\alpha$,
\[
\tilde{\T}(\alpha) = \int \tilde{f}(\beta, \alpha) \, \tilde{R}(d\beta \given \alpha) = 2\int \beta \, \tilde{R}(d\beta \given \alpha) + 2\alpha + 2c
\]
is greater than or equal to $c$. Thus the kernel
\begin{equation}
\tilde{T}(d\beta \given \alpha) = \frac{(\beta + \alpha  + c) \, \tilde{R}(d\beta \given \alpha)}{\int (v+ \alpha  + c) \, \tilde{R}(dv \given \alpha)} \label{love.kernel}
\end{equation}
defines a $\tilde{\mu}$-symmetric Markov chain on $\nnr$ where
\[ 
\tilde{\mu}(A) = \int_{A} \tilde{\T}(\alpha) \, \tilde{\nu}(d\alpha)
\; \qquad A \in \A \; .
\] 
Clearly, for any integer $n$, $\tilde{\mu}([0,n)) < \infty$.

The next step is to verify the conditions of Theorem~\ref{thm:recurrence} which will imply the Markov chain associated with $\tilde{T}$ is locally $\tilde{\mu}$-recurrent.  Since $\tilde{\mu}$ and $\tilde{\nu}$ are equivalent measures we will also conclude that the chain is locally $\tilde{\nu}$-recurrent.  In Appendix \ref{app:drift.calculations} it is shown that  $\tilde{T}([0,m] \given \alpha)$ is continuous as a function of $\alpha$ implying condition \eqref{sup.T} of Theorem \ref{thm:recurrence}. For $k$ a nonnegative integer set
\begin{equation*}
m_{k}(\alpha) = \int (\beta-\alpha)^{k} \tilde{T} (d\beta \given \alpha) \; .
\end{equation*}
Additional calculations given in Appendix \ref{app:drift.calculations} show that
\begin{equation} \label{eq:first.moment}
\frac{m_{3}(\alpha)}{m_{2}(\alpha )} = O(1) \quad \text{as}~\alpha \to \infty 
\end{equation}
and
\begin{equation}
m_{1}(\alpha) = \frac{8\alpha + \psi_{1}(\alpha)}{\int (\beta + \alpha + c) \, \tilde{R}(d\beta \vert \alpha)} \label{drift.lhs}
\end{equation}
where $\psi_{1}(\alpha)=O(1)$ as $\alpha \to \infty$.  Letting $\phi(\alpha) = 1/\sqrt{\alpha}$, note that $\lim_{\alpha \to \infty}\phi(\alpha)\log(\alpha)=0$ and, by calculations in Appendix~\ref{app:drift.calculations},
\begin{equation}
\frac{m_{2}(\alpha)}{2\alpha} [1 + \phi(\alpha)]   = \frac{8\alpha + 8\sqrt{\alpha}+\psi_{2}(\alpha)}{\int (\beta + \alpha + c) \, \tilde{R}(d\beta \vert \alpha)} \label{drift.rhs}
\end{equation}
where $\psi_{2}(\alpha)=O(1)$ as $\alpha \to \infty$.  It is clear by inspection of \eqref{eq:first.moment}, \eqref{drift.lhs}, and \eqref{drift.rhs} that the conditions \eqref{moment.cond.1} and \eqref{moment.cond.2} of Theorem~\ref{thm:recurrence} are satisfied for $\alpha$ large enough.  Hence the chain is locally $\tilde{\nu}$-recurrent.

By Theorem \ref{thm:reduced result}, the family of priors $\Fnu$ is
strongly admissible with respect to $\Phi_{f}$.

\end{proof}

\appendix

\section{Preliminaries}

We begin by stating some existing results concerning local recurrence
and introduce Blyth's method.  This material plays a fundamental role
in our proofs of Theorems~\ref{thm:main result},~\ref{thm:family}
and~\ref{thm:reduced result}.

\subsection{Local Recurrence}
The purpose here is to give two characterizations of local recurrence
for general symmetric Markov chains, hence the notation is consistent
with that of Sections~\ref{sec:preliminaries} and~\ref{sec:nnr}.

Let $\statespace$ be a Polish space and denote the Borel
$\sigma$-algebra by $\D$.  Let $K : \D \times \statespace \to [0,1]$
be a Markov kernel on $(\statespace,\D)$. Let $\xi$ be a non-trivial
$\sigma$-finite measure and recall that a set is $\xi$-proper if its
measure under $\xi$ is positive and finite.  Throughout this section
$K$ is assumed to be $\xi$-symmetric.

\begin{theorem}[\citeauthor{Eaton:2004} \cite{Eaton:2004}] \label{increasing.Cs}
The Markov chain $W$ is locally $\xi$-recurrent if and only if there
exists a sequence of $\xi$-proper sets increasing to the state space
such that each is locally $\xi$-recurrent. 
\end{theorem}

Let $L^{2}(\xi)$ be the space of $\xi$-square integrable
functions. Then the quantity
\begin{equation}
\label{eq:Dirichlet form}
\Delta (h; K, \xi) = \frac{1}{2} \iint (h(\theta) - h(\eta))^{2}
K(d\theta | \eta) \xi(d\eta) \hspace*{9mm} h \in L^{2}(\xi) 
\end{equation}
is called a \textit{Dirichlet form}.  Also, if $D$ is a $\xi$-proper set, define
\[ 
\mathcal{H}_{\xi}(D) := \{ h \in L^{2}(\xi) :  h \geq I_{D}\}  
\]
where $I_{D}$ is the indicator function of the set $D$.  A
characterization of local $\xi$-recurrence in terms of $\Delta$ is
given by the following result. 

\begin{theorem}[\citeauthor{Eaton:2001} \cite{Eaton:2001}]
\label{joeppendix}
A set $D \in \D$ is locally $\xi$-recurrent if and only if
\[
\inf_{\mathcal{H}_{\xi}(D)} \Delta(h;K,\xi) = 0 \; .
\]
\end{theorem}

\subsection{Blyth's Method}
\label{app:connection}
Consider the posterior distributions obtained from perturbations of
the prior measure $\nu$.  Let $g : \parmspace \to \mathbb{R}^{+}$ be
such that the perturbation
\begin{equation}
\label{eq:perturbed prior}
\nu_{g}(d\theta) = g(\theta) \nu(d\theta)
\end{equation}
is $\sigma$-finite.  Assume
\begin{equation}
\label{eq:ghat}
\hat{g}(x) := \int g(\theta) Q(d\theta | x) < \infty \; .
\end{equation}
Letting
\[
M_{g}(dx) := \int P(dx|\theta) g(\theta) \nu(d\theta) 
\]
it is easy to see that $M_{g}(dx)= \hat{g}(x) M(dx)$ is
$\sigma$-finite and that the posterior obtained from the perturbed
prior is 
\begin{equation}
\label{eq:perturbed posterior}
Q_{g}(d\theta | x) = \frac{g(\theta)}{\hat{g}(x)} Q(d\theta | x) \; .
\end{equation}
Let $L^{1}(\nu)$ be the set of all $\nu$-integrable functions and define
\begin{equation}
\label{eq:Gnu}
\kernels = \{ g \in L^{1}(\nu) : g \geq 0, \; g\; \text{bounded}
\text{ and } \nu_{g} (\parmspace) > 0  \} \; .  
\end{equation}
Let $I_{C}$ be the indicator function of the set $C$.  Of particular
interest are the subfamilies  
\[
\kernels(C) = \{ g \in \kernels : g \geq I_{C} \} \; .
\]
Let $\hat{\varphi}$ be the Bayes estimator of $\varphi$ under the
prior $\nu$, and let $\hat{\varphi}_{g}$ be the Bayes estimator of
$\varphi$ under the prior $\nu_{g}$.  Also recall the definition of
the risk function $\risk$ at \eqref{eq:risk}.  A key quantity in
connecting Markov chains to admissibility is the \textit{integrated
  risk difference} at $g \in \kernels$ with respect to $\nu$ against
$\hat{\phee}$:
\begin{equation}
\label{eq:IRD}
\ird_{\nu}(g;\hat{\phee}) := \int \left[ \risk(\theta;\hat{\varphi}) -
  \risk(\theta;\hat{\varphi}_{g}) \right] \,\nu_{g}(d\theta).  
\end{equation}
Notice that the integrability assumptions on $\risk$ in
Assumption~\ref{ass:risk} ensure that $\ird_{\nu}$ is
well-defined. The connection of $\ird_{\nu}$ with almost-$\nu$
admissibility is given by Blyth's method.
\begin{theorem}[Blyth's Method] \label{thm:Blyth}
Let $\hat{\phee}$ be an estimator and $\nu$ a $\sigma$-finite measure
on the parameter space. If  $\inf_{\kernels(C)}
\ird_{\nu}(g;\hat{\phee}) = 0$ for every  $C \subseteq \parmspace$
such that $0 < \nu(C) <\infty$, then $\hat{\phee}$ is almost-$\nu$
admissible. 
\end{theorem}

\section{Proof of Theorem~\ref{thm:main result}} 
\label{app:main result}
We will develop a connection between the kernel $T$ defined in
\eqref{T.def} and the $\ird_{\nu}$ which will be key to proving
Theorem~\ref{thm:main result} via Theorems~\ref{joeppendix}
and~\ref{thm:Blyth}. We begin with some preliminary results before we
prove Theorem~\ref{thm:main result}. 

\subsection{Preliminary Results}\label{app:prelims}

Recall $\mathcal{T}(\eta)=\psi(\parmspace | \eta)$ where $\psi$ is
defined at \eqref{psi.def}.  Define  
\begin{equation}
 \mu(C) = \int_{C} \mathcal{T}(\eta) \, \nu(d\eta) \quad \text{for
   all}~C \in \C  \; . \label{mu.def} 
\end{equation}

\begin{proposition} \label{rnd.thm.1} Suppose Assumption~\ref{ass:T}  and condition \eqref{density.symmetry} hold. Then
\begin{enumerate}
\item[(a)] the measure $\mu$ is $\sigma$-finite and equivalent to $\nu$,
\item[(b)] $\mathcal{T}$ is a Radon-Nikodym derivative of $\mu$ with respect to $\nu$, 
\item[(c)] $T$ is a $\mu$-symmetric Markov transition kernel,
\item[(d)] $T(d\theta | \eta) \mu(d\eta) = f(\theta, \eta) R(d\theta | \eta) \nu(d\eta)$,
\item[(e)] if $\mu'$ is proportional to $\mu$, $T$ is a $\mu'$-symmetric Markov transition kernel.
\end{enumerate}
\end{proposition}

\begin{proof}
That $\mu$ is $\sigma$-finite follows easily from Assumption~\ref{ass:T}.  Since $\mathcal{T}$ is a nonnegative measurable function, for any $\nu$-null set $C$,
\[ 
\mu(C) = \int_{C} \mathcal{T}(\eta) \, \nu(d\eta) = 0 \; .
\]
Furthermore, since $\T(\eta) > 0$ for all $\eta$ in the parameter space, every $\nu$-positive set is a $\mu$-positive set.  Thus, if $C$ is $\mu$-null, it is also $\nu$-null.  Therefore, $\mu$ and $\nu$ are equivalent measures.  Furthermore, $\T$ is a Radon-Nikodym derivative of $\mu$ with respect to $\nu$ since $\mu$ is absolutely continuous with respect to $\nu$, and by assumption $\T$ is a nonnegative measurable function such that for any measurable set $C$, $\mu(C)$ is given by equation \eqref{mu.def}.

At any point $\eta$ in the parameter space, since $\psi( \cdot \given
\eta)$ is a nontrivial finite measure, normalizing by
$\T(\eta)=\psi(\parmspace \given \eta)$ produces a probability
measure. Recall that $\psi(C \given \cdot)$ is a nonnegative
measurable function where $C$ is an element of the Borel sets $\C$.
Since $\parmspace$ is a Borel measurable set, $\T$ is a Borel
measurable function.  Also by hypothesis, $\T (\eta)$ is positive and
finite for all $\eta$.  The reciprocal function is continuous, hence
Borel measurable, on the positive real numbers, and the composition $1
/ \T$ of Borel measurable functions is a Borel measurable function.
Thus, the product of Borel measurable functions $1 / \T$ and $\psi(C
\given \cdot)$ is a Borel measurable function for any $C \in \C$.   

Let $A$ and $B$ be measurable sets.  Since
\begin{align*}
S(A,B)
&= \iint I_{A}(\eta) I_{B}(\theta) \, T(d\theta \given \eta) \,
\mu(d\eta) \\ 
&= \iint I_{A}(\eta) I_{B}(\theta) \, f(\theta, \eta) \, R(d\theta
\given \eta) \, \nu(d\eta) && \text{by \eqref{psi.def}, \eqref{T.def},
  and \eqref{mu.def}} \\ 
&= \iint I_{B}(\theta) I_{A}(\eta) \, f(\eta, \theta) \, R(d\eta
\given \theta) \, \nu(d\theta) && \text{by \eqref{detailed.balance}
  and \eqref{density.symmetry}} \\ 
&= \iint I_{B}(\theta) I_{A}(\eta) \, T(d\eta \given \theta) \,
\mu(d\theta) && \text{by \eqref{psi.def}, \eqref{T.def}, and
  \eqref{mu.def}} \\ 
&= S(B,A),
\end{align*}
the kernel $T$ is $\mu$-symmetric.
Note that the substitution
\begin{equation}
T(d\theta \given \eta)  \, \mu(d\eta) = f(\theta, \eta) \, R(d\theta
\given \eta) \, \nu(d\eta) \label{sub} 
\end{equation}
holds as a consequence.  

Let $\mu' = c \mu$ for some $c > 0$.  Then by \eqref{detailed.balance}
and \eqref{sub} 
\begin{align*}
T(d\theta | \eta) \mu'(d\eta)
& = cT(d\theta|\eta) \mu(d\eta) \\
& = c f(\theta, \eta) R(d\theta | \eta) \nu(d\eta)\\
& = c f(\eta, \theta) R(d\eta | \theta) \nu(d\theta) \\
& = c T(d\eta | \theta) \mu(d\theta) \\
& = T(d\eta | \theta) \mu'(d\theta)
\end{align*}
which implies $T$ is $\mu'$-symmetric.  
\end{proof}

We can now develop a connection between the Markov kernel $T$ and the
integrated risk difference $\ird_{\nu}$. Our argument will require the
following known result; recall the definition of $\kernels$ from
\eqref{eq:Gnu}. 

\begin{proposition}[\citeauthor{Eaton:2001}, \citeyear{Eaton:2001}] 
\label{prop:Link}
If $g \in \mathcal{G}_{\nu}$ and $R$ is Eaton's kernel, then
\begin{equation*}
\ird_{\nu}(g;\hat{\varphi}) \leq \iint \lVert
\varphi(\theta)-\varphi(\eta)\rVert^{2}
(\sqrt{g}(\theta)-\sqrt{g}(\eta))^{2} \, R(d\theta \given \eta) \,
\nu(d\eta)\; .  
\end{equation*}
\end{proposition}

\noindent By Proposition~\ref{rnd.thm.1} we see that $T$ is
$\mu$-symmetric so that the relevant Dirichlet form, recall
\eqref{eq:Dirichlet form}, is
\[ 
\Delta(h; T,\mu) = \frac{1}{2} \iint ( h(\eta) - h(\theta) )^{2} \,
T(d\theta \given \eta) \, \mu(d\eta)  \qquad h \in L^{2}(\mu) \; . 
\]

\begin{proposition} \label{prop:irdledf} Suppose
  Assumption~\ref{ass:T} and condition \eqref{density.symmetry}
  hold. If $\varphi \in \varPhi_{f}$ and $g \in \mathcal{G}_{\nu}$,
  then there is a measure $\mu_{\varphi}$ which is proportional to
  $\mu$ such that
\[ 
\ird_{\nu} (g;\hat{\varphi}) \leq 2 \Delta(\sqrt{g}; T, \mu_{\varphi}) \; .
\]
\end{proposition}

\begin{proof}
Since $\varphi \in \varPhi_{f}$ there exists $0 < M_{\varphi} < \infty$ such that
\[
\|\varphi(\theta) - \varphi(\eta)\|^{2} \le M_{\varphi} f(\theta, \eta) \quad \text{ for all } \theta, \eta
\]
For $C \in \C$ define $\mu_{\varphi}(C) = M_{\varphi} \mu(C)$ and
suppose $g \in \mathcal{G}_{\nu}$.  In the following the first
inequality is from Proposition~\ref{prop:Link} and the second is
obtained from Proposition~\ref{rnd.thm.1}(d) and that $\varphi \in
\varPhi_{f}$
\begin{align*}
 \ird_{\nu} (g;\hat{\varphi}) 
&\leq \iint \lVert\phee(\theta)-\phee(\eta)\rVert^{2} (\sqrt{g}(\theta)-\sqrt{g}(\eta))^{2} \, R(d\theta \given \eta) \, \nu(d\eta) \\
&\leq  \iint (\sqrt{g}(\theta)-\sqrt{g}(\eta))^{2} \, T(d\theta \given \eta)  \, M_{\varphi} \mu(d\eta) \\ 
& = 2 \Delta(\sqrt{g}; T, \mu_{\varphi})  \; .
\end{align*}
\end{proof}

\subsection{Proof of Theorem~\ref{thm:main result}}
Let the measures $\psi( \cdot \vert \eta)$ and $\mu$ and the transition kernel $T$ be as defined at equations \eqref{psi.def}, \eqref{mu.def}, and \eqref{T.def}, respectively. Also, suppose $\varphi \in \varPhi_{f}$.  By Proposition~\ref{prop:irdledf} if $g \in \mathcal{G}_{\nu}$, then
\[ 
\ird_{\nu} (g;\hat{\varphi}) \leq 2 \Delta(\sqrt{g}; T, \mu_{\varphi}) 
\]
where $\mu_{\varphi}= M_{\varphi} \mu$.   By Proposition~\ref{rnd.thm.1} we have that $\mu$ and $\nu$ are equivalent measures and $T$ is a $\mu_{\varphi}$-symmetric Markov kernel.  Since the measures $\mu_{\varphi}$ and $\nu$ are equivalent and the chain with kernel $T$ is locally $\nu$-recurrent, it is also locally $\mu_{\varphi}$-recurrent---that is, every $\mu_{\varphi}$-proper set $C$ is locally $\mu_{\varphi}$-recurrent. Thus, by Theorem~\ref{joeppendix}, 
\[
\inf_{\mathcal{H}_{\mu_{\varphi}}(C)} \Delta(h;T,\mu_{\varphi}) = 0
\]
where $\mathcal{H}_{\mu_{\varphi}}(C)$ collects the square-integrable dominators of $I_{C}$.

%The set $\mathcal{G}_{\nu_{u}}$ of distribution kernels with respect to the %perturbed prior is the same as $\kernels$ since for all $g \in L^{1}(\nu)$
%\[ 
%\frac{1}{c} \int \lvert g (\theta) \rvert \, \nu(d\theta)  \leq \int \lvert g(\theta) %\rvert \, \nu_{u}(d\theta) \leq  c \int \lvert g(\theta) \rvert \, \nu(d\theta) \; .
%\]
Now let $h \in \mathcal{H}_{\mu_{\varphi}}(C)$ and recall that by assumption $\T$ is a measurable function uniformly bounded away from zero $\nu$-almost everywhere. Hence there exists some $\epsilon > 0$ such that 
\begin{align*}
\epsilon \int h^{2}(\eta) \, \nu (d\eta)  & \le \int h^{2} (\eta) \T(\eta) \, \nu(d\eta) \\ 
& = \int h^{2}(\eta) \, \mu (d\eta) \\
& = \frac{1}{M_{\varphi}} \int h^{2}(\eta) \, \mu_{\varphi} (d\eta) \\
&< \infty \, .
\end{align*}
Hence $h \in L^{2}(\nu)$ and we conclude that $\mathcal{H}_{\mu_{\varphi}}(C) \subseteq \mathcal{G}_{\nu}(C)$.  Moreover, if $\sqrt{h} \in \mathcal{H}_{\mu_{\varphi}}(C)$, then $h \in \mathcal{G}_{\nu}(C)$. Thus we obtain
\[
\inf_{\mathcal{G}_{\nu}(C)} \ird_{\nu}(g;\hat{\phee}) = 0 .
\]
Therefore, $\hat{\phee}$ is an almost-$\nu$ admissible estimator by Theorem~\ref{thm:Blyth}.  

\subsection{Proof of Theorem~\ref{thm:family}}

Recall that $\nu_{u} \in \Fnu$ is a bounded perturbation of $\nu$. The sampling and posterior distributions define the $\nu_{u}$-symmetric Eaton kernel 
\[ 
R_{u} (d\theta \given \eta) = \int_{\X} Q_{u}(d\theta \given x) \, P(dx \given \eta) 
\]
where the perturbed posterior was defined at \eqref{eq:perturbed posterior}.  The mean of $\phee$ with respect to $Q_{u}( \cdot \given x)$ is the formal Bayes estimator of $\phee(\theta)$ under squared error loss. We denote it $\hat{\phee}_{u}$ to emphasize its dependence on the perturbed prior.  

Since $\nu_{u} \in \Fnu$, there exists $0 < c < \infty$ such that $\nu_{u}(d\theta) = u(\theta) \nu(d\theta)$ and $1/c < u(\theta) < c$.  Recall the definition of $\mu_{\phee}$ from the proof of Proposition~\ref{prop:irdledf}.   In the following the first inequality is from Proposition~\ref{prop:Link}, the second follows by noting that $u$ and $1/\hat{u}$ are both bounded above by $c$ while the third is obtained from Proposition~\ref{rnd.thm.1}(d) and that $\varphi \in \varPhi_{f}$
\begin{align*}
 \ird_{\nu_{u}} (g;\hat{\varphi}_{u}) 
&\leq \iint \lVert\phee(\theta)-\phee(\eta)\rVert^{2} (\sqrt{g}(\theta)-\sqrt{g}(\eta))^{2} \, R_{u}(d\theta \given \eta) \, \nu_{u}(d\eta) \\
& \le c^{3} \iint \lVert\phee(\theta)-\phee(\eta)\rVert^{2} (\sqrt{g}(\theta)-\sqrt{g}(\eta))^{2} \, R(d\theta \given \eta) \, \nu(d\eta) \\
&\leq c^{3} \iint (\sqrt{g}(\theta)-\sqrt{g}(\eta))^{2} \, T(d\theta \given \eta)  \, M_{\varphi} \mu(d\eta) \\ 
& = 2 c^{3} \Delta(\sqrt{g}; T, \mu_{\varphi})  \; .
\end{align*}
The remainder of the proof follows the proof of Theorem~\ref{thm:main result} exactly.

\section{Proof of Theorem~\ref{thm:reduced result}}
\label{app:reduced result}
Before proving Theorem~\ref{thm:reduced result} we require an analogue
of Proposition~\ref{rnd.thm.1}.  Given a nonnegative function
$\tilde{h}$ on $\nnr$, define
\begin{equation}
h(\theta) = \tilde{h}(t(\theta)). \label{h.tilde}
\end{equation}

\begin{proposition} \label{rnd.thm.2}
Suppose Assumptions~\ref{ass:partition} and~\ref{ass:tildeT} hold.  If $\tilde{f}(\beta, \alpha) = \tilde{f}(\alpha, \beta)$, then
\begin{enumerate}
\item[(a)] $\tilde{\mu}$ is $\sigma$-finite and equivalent to $\tilde{\nu}$,
\item[(b)] $\tilde{\T}$ is a Radon-Nikodym derivative of $\tilde{\mu}$ with respect to $\tilde{\nu}$, 
\item[(c)] $\tilde{T}$ is a $\tilde{\mu}$-symmetric Markov transition kernel,
\item[(d)] $\tilde{T}(d\beta \given \alpha) \, \tilde{\mu}(d\alpha) = \tilde{f}(\beta, \alpha) \, \tilde{R}(d\beta \given \alpha) \, \tilde{\nu}(d\alpha)$,
\item[(e)] $\tilde{h} \in L^{2}(\tilde{\mu})$ implies $h$ is in $L^{2}(\mu)$, \item[(f)]$\Delta(\tilde{h} ;\tilde{T},\tilde{\mu}) = \Delta(h ;T,\mu)$, and
\item[(g)] $\mathcal{T}(\eta)$ is almost-$\nu$ bounded away from zero.
\end{enumerate}
\end{proposition}

\begin{proof}
The proof of the first 3 assertions follows exactly the proof of the first 3 assertions in Proposition~\ref{rnd.thm.1} with $\nnr$, $\A$, $\alpha$, $\beta$, $\tilde{\nu}$, and $\tilde{R}$ substituted for $\parmspace$, $\C$, $\eta$, $\theta$, $\nu$, and $R$, respectively. 
The substitution
\begin{equation}
\tilde{T}(d\beta \given \alpha) \, \tilde{\mu}(d\alpha) = \tilde{f}(\beta, \alpha) \, \tilde{R}(d\beta \given \alpha) \, \tilde{\nu}(d\alpha) \label{reduced.sub}
\end{equation}
follows as a consequence.
If $\tilde{h} \in L^{2}(\tilde{\mu})$, then $h \in L^{2}(\mu)$ since 
\begin{align*}
\int \tilde{h}^{2}(\alpha) \tilde{\mu}(d\alpha)
&= \int \tilde{h}^{2}(\alpha) \tilde{\T}(\alpha) \tilde{\nu}(d\alpha) \\
&= \int \tilde{h}^{2}(\alpha) \int \tilde{f}(\beta, \alpha) \int \tilde{Q}(d\beta \given x) \tilde{P}(dx \given \alpha) \tilde{\nu}(d\alpha) \\
&= \int h^{2}(\eta) \int f(t(\theta), t(\eta)) \int Q(d\theta \given x) P(dx \given \eta) \nu(d\eta) \\
&= \int h^{2}(\eta) \T(\eta) \nu(d\eta) \\
&= \int h^{2}(\eta) \mu(d\eta). 
\end{align*}
Since
\begin{align*}
&\quad \iint (\tilde{h}(\beta) - \tilde{h}(\alpha))^{2} \tilde{T}(d\beta \given \alpha) \tilde{\mu}(d\alpha)\\
&= \iint (\tilde{h}(\beta) - \tilde{h}(\alpha))^{2} \tilde{f}(\beta, \alpha)  \tilde{R}(d\beta \given \alpha)  \tilde{\nu}(d\alpha) && \text{by \eqref{reduced.sub}} \\
&= \iint (\tilde{h}(\beta) - \tilde{h}(\alpha))^{2} \tilde{f}(\beta, \alpha)  \int \tilde{Q}(d\beta \given x) \tilde{P}(dx \given \alpha) \tilde{\nu}(d\alpha) && \text{by \eqref{reduced.Eaton.kernel}} \\
&= \iint (h(\theta) - h(\eta))^{2} f(t(\theta), t(\eta)) \int
Q(d\theta \given x) P(dx \given \eta) \nu(d\eta)  &&\text{by
  \eqref{nu.tilde}, \eqref{p.tilde}, \eqref{q.tilde}, and
  \eqref{h.tilde}} \\
&= \iint (h(\theta) - h(\eta))^{2} f(t(\theta), t(\eta)) R(d\theta \given \eta) \nu(d\eta) && \text{by \eqref{Eaton.kernel}} \\
&= \iint (h(\theta) - h(\eta))^{2} T(d\theta \given \eta) \mu(d\eta) && \text{by \eqref{sub}},
\end{align*}
it follows that $\Delta(\tilde{h} ;\tilde{T},\tilde{\mu}) = \Delta(h ;T,\mu)$.  Finally, recall that $\tilde{\T}$ is almost-$\tilde{\nu}$ bounded away from zero by some positive constant $c$. 
Note that for any $\tilde{\nu}$-proper set $A$,
\[ 
\int_{t^{-1}(A)} \T(\eta) \nu(d\eta) = \int_{A} \tilde{\T}(\alpha) \tilde{\nu}(d\alpha) \geq c \tilde{\nu}(A) = c \nu(t^{-1}(A)) . 
\]
It follows that $\T$ is almost-$\nu$ bounded away from zero.
\end{proof}

We are now ready to prove Theorem~\ref{thm:reduced result}.

\subsection{Proof of Theorem~\ref{thm:reduced result}}
By Proposition~\ref{rnd.thm.2}, the measures $\tilde{\mu}$ and $\tilde{\nu}$ are equivalent, and the kernel $\tilde{T}$ is $\tilde{\mu}$-symmetric.
Thus, the chain with kernel $\tilde{T}$ is locally $\tilde{\mu}$-recurrent, and by Theorem~\ref{joeppendix}
\[
\inf_{\mathcal{H}_{\tilde{\mu}}(A)} \Delta(\tilde{h};\tilde{T},\tilde{\mu}) = 0
\]
for any $\tilde{\mu}$-proper set $A$.

Let $(A_{i})$ be a sequence of $\tilde{\mu}$-proper sets increasing to $\nnr$.
Letting $C_{i} = t^{-1}(A_{i})$ defines a sequence of $\mu$-proper sets increasing to $\parmspace$.
Since $h(\theta)=\tilde{h}(t(\theta))$, $\tilde{h} \geq I_{A_{i}}$ implies that $h \geq I_{C_{i}}$.
By Proposition~\ref{rnd.thm.2}, $\tilde{h} \in L^{2}(\tilde{\mu})$ implies $h \in  L^{2}(\mu)$ and the corresponding Dirichlet forms are equal.
Thus, 
\[
\inf_{\mathcal{H}_{\mu}(C_{i})} \Delta(h;T,\mu) = \inf_{\mathcal{H}_{\tilde{\mu}}(A_{i})} \Delta(\tilde{h};\tilde{T},\tilde{\mu}) = 0
\]
 for every $C_{i}$ in the sequence.
 By Theorem~\ref{increasing.Cs}, the chain with kernel $T$ is locally $\mu$-recurrent.
 Since $\mu$ and $\nu$ are equivalent measures, the chain is locally $\nu$-recurrent as well.

 By Proposition~\ref{rnd.thm.2}, $\T$ is almost-$\nu$ bounded away from zero, and thus all of the conditions for Theorem~\ref{thm:family} are satisfied.  Therefore, under squared error loss, the family $\Fnu$ is $\varPhi_{f}$-admissible.

\section{Multivariate Normal}\label{app:drift.calculations}

Here we give the supplemental arguments required for the proof of Theorem~\ref{illustration}. 

\subsection{Existence of the integrated risk difference} In order to appeal to Blyth's method, we need to know that the integrated risk differences are defined.
Otherwise the bounding inequality \eqref{prop:Link} for the integrated risk differences is meaningless.
The proper Bayes estimators necessarily have finite integrated risks, so it is sufficient to show that 
\[ \int \risk(\theta;\hat{\theta}) \, g_{n}(\theta) \, \nu(d\theta) < \infty \]
where $g_{n}$ is the indicator of $C_{n}$ a closed ball around zero with radius $n$.
Since $\nu$ is $\sigma$-finite, $\nu(C_{n})$ is finite.
Since the risk function is real-valued and continuous, it attains a finite maximum on $C_{n}$; call it $M_{n}$.
Therefore, 
\[ \int  \risk(\theta;\hat{\theta}) \, g_{n}(\theta) \, \nu(d\theta) \leq M_{n} \nu(C_{n}) < \infty \]
and the integrated risk differences are defined.

\subsection{Continuity of the transition kernel} 
Let $m$ be a nonnegative integer.
We wish to show that $T([0,m] \given \alpha)$ is continuous as a function of $\alpha$.
Since by definition
\[ T([0,m] \given \alpha) = \frac{\int_{0}^{m} \beta \tilde{R}(d\beta \given \alpha) + \alpha \tilde{R}([0,m] \given \alpha) + c \tilde{R}([0,m] \given \alpha)}{\int \beta \tilde{R}(d\beta \given \alpha) + \alpha + c}  \]
it is sufficient to show that $\tilde{R}([0,m] \given \alpha)$ is continuous and that continuity of 
\[ \int_{0}^{m} \beta \tilde{R}(d\beta \given \alpha) \quad \text{and} \quad \int \beta \tilde{R}(d\beta \given \alpha) \]
follows from there.

Fix $\alpha^{*} \in [0,\infty)$ and $\delta$ greater than zero.
Let $S_{\delta}(\alpha^{*}) = [0,(\sqrt{\alpha^{*}}+ \delta)^{2})$.
Let $\left( \alpha_{n} \right)$ be a sequence with limit $\alpha^{*}$ whose elements are in $S_{\delta}(\alpha^{*})$.
Define
\[
A_{\delta}(\alpha^{*}) = \{ x : \lVert x \rVert < 2(\sqrt{\alpha^{*}}+\delta) \} \,
\]
---an open ball in $\X$ with radius greater than twice $\sqrt{\alpha^{*}}$.
Let $\Xi$ denote the surface of the unit hypersphere in $\mathbb{R}^{p}$.
For any $x \in A_{\delta}^{c}(\alpha^{*})$---that is, any point such that $\lVert x \rVert$ is greater than twice $\sqrt{\alpha^{*}}+\delta$---and $\xi$ on the unit hypersphere
\[
\lVert x - \xi\sqrt{\alpha_{n}} \rVert \geq \left\lVert x \left( 1 -  \sqrt{\alpha_{n}} / \lVert x \rVert \right) \right\rVert \geq \left\lVert x/2  \right\rVert
\]
since $x$ cannot be closer to any point with radius $\sqrt{\alpha_{n}}$ than it is to $x\sqrt{\alpha_{n}}/\lVert x \rVert$, the point with radius $\sqrt{\alpha_{n}}$ on the common ray and since $\sqrt{\alpha_{n}}$ is less than $\sqrt{\alpha^{*}}+\delta$ by construction.
Recall that 
\[
\tilde{p}(x \vert \alpha) = \int_{\Xi} (2\pi)^{-p/2} e^{-\frac{1}{2} \lVert x - \xi\sqrt{\alpha} \rVert} \, \pi(d\xi)
\]
where $\pi$ is the uniform distribution on $\Xi$.
Let 
\[
g_{1}(x) = (2\pi)^{-p/2} \left[ I_{A_{\delta}(\alpha^{*})}(x) + I_{A^{c}_{\delta}(\alpha^{*})}(x) e^{-\frac{1}{8} \lVert  x \rVert^{2} } .\right]
\]
Since $\int g_{1}(x) \, \pi(d\xi) < \infty$ and $g_{1}(x)$ dominates the integrand of $\tilde{p}(x \vert \alpha_{n})$, we can say that
\[ \lim_{n \rightarrow \infty} \tilde{p}(x \vert \alpha_{n}) = \tilde{p}(x \vert \alpha^{*}) \, . \]
Furthermore, $\int g_{1}(x) \, dx < \infty$ and $g_{1}(x)$ dominates $\tilde{p}(x \vert \alpha_{n})$ by monotonicity of the integral.
Therefore,
\begin{align*}
\lim_{n \rightarrow \infty}  \int \tilde{p}(x \vert \alpha_{n}) \, dx = \int \lim_{n \rightarrow \infty} \tilde{p}(x \vert \alpha_{n}) \, dx = \int \tilde{p}(x \vert \alpha^{*}) \, dx \, .
\end{align*}

Denote the Borel $\sigma$-algebra on $[0,\infty)$ by $\mathcal{B}$, and choose $B \in \mathcal{B}$.
Since the densities $\tilde{p}$ and $\tilde{q}$ are necessarily nonnegative and measurable, Fubini says
\[
\tilde{R}(B \vert \alpha) = \int_{B} \int_{\mathcal{X}} \tilde{Q}(d\beta \vert x) \, \tilde{P}(dx \vert \alpha) = \int_{\mathcal{X}} \tilde{p}(x \vert \alpha) \, \tilde{Q}(B \vert x) \, dx \, .
\]
Similarly,
\[
\int_{B} \beta \tilde{R}(d\beta \vert \alpha) = \int_{B} \beta \, \int_{\mathcal{X}} \tilde{Q}(d\beta \vert x) \, \tilde{P}(dx \vert \alpha) = \int_{\mathcal{X}} \tilde{p}(x \vert \alpha) \int_{B} \beta  \, \tilde{q}(\beta \vert x) \, d\beta \, dx \, .
\]
Let
\[
f(x) = \tilde{p}(x \vert \alpha^{*}) \tilde{Q}(B \vert x) \quad \text{and} \quad f_{n}(x) = \tilde{p}(x \vert \alpha_{n}) \tilde{Q}(B \vert x)
\]
so that 
\[ \lim_{n \rightarrow \infty} f_{n}(x) = f(x). \]
Since $\tilde{Q}(B \vert x)$ is a probability,
\[
f_{n}(x) \leq \tilde{p}(x \vert \alpha_{n}) \leq g_{1}(x) \, .
\]
Hence, by the dominated convergence theorem,
\[
\lim_{n \rightarrow \infty} \int f_{n}(x) \, dx = \int \lim_{n \rightarrow \infty} f_{n}(x) \, dx = \int  f(x) \, dx .
\]
That is,
\[  
\lim_{n \rightarrow \infty} \tilde{R}(B \vert \alpha_{n}) = \tilde{R}(B \vert \alpha^{*}).
\]
Therefore, $\tilde{R}(B \vert \alpha)$ is a continuous function of $\alpha$.

Now note that
\[
\int_{B} \beta \, \tilde{q}(\beta \vert x) \, d\beta \leq \int \beta \, \tilde{q}(\beta \vert x) \, d\beta \, .
\]
By Proposition A.7 of \citet{Eaton:2008}, for some bounded $\psi(\lVert x \rVert^{2})$,
\[
\int \beta \, \tilde{q}(\beta \vert x) \, d\beta = -p + \lVert x \rVert^{2} + \psi(\lVert x \rVert^{2})
\]
so that there exists a constant $k >0$ such that
\[
\int \beta \, \tilde{q}(\beta \vert x) \, d\beta \leq k + \lVert x \rVert^{2} \, .
\]
Let $g_{2}(x) = k + \lVert x \rVert^{2}$, and let $g(x) = g_{1}(x)g_{2}(x)$.
Let 
\[
f_{n}(x) = \tilde{p}(x \vert \alpha_{n}) \int_{B} \beta  \, \tilde{q}(\beta \vert x) \, d\beta \quad \text{and} \quad f(x) = \tilde{p}(x \vert \alpha^{*}) \int_{B} \beta  \, \tilde{q}(\beta \vert x) \, d\beta
\]
Note that 
\[
\lim_{n \rightarrow \infty} f_{n}(x) = f(x)
\] 
and $g(x)$ dominates $f_{n}(x)$.
One can show that
\[
\int g(x) \, dx < \infty 
\]
by expanding the product $g_{1}(x)g_{2}(x)$ and integrating the components.
Now, by the Dominated Convergence Theorem,
\[
\lim_{n \rightarrow \infty} \int f_{n}(x) \, dx = \int \lim_{n \rightarrow \infty} f_{n}(x) \, dx = \int f(x) \, dx \, .
\]
That is, 
\[
\lim_{n \rightarrow \infty} \int_{B}\beta R(d\beta\vert\alpha_{n}) = \int_{B}\beta R(d\beta\vert\alpha^{*}) \, .
\]
Therefore, $\int_{B}\beta R(d\beta\vert\alpha)$ is continuous as a function of $\alpha$.
Since $B$ is any Borel set, the continuity holds for both $[0,m]$ and $\nnr$.
Finally, since $T([0,m] \given \alpha)$ is an algebraic combination of continuous functions, it is itself a continuous function. 

\subsection{Moment conditions}

We are interested in moments of the transition kernel about the current state:
\[
m_{k}(\alpha) =  \int (\beta - \alpha)^{k} \, \tilde{T}(d\beta \vert \alpha) \, .
\]
We can express these moments in terms of $\tilde{R}$ as
\[
m_{k}(\alpha) = \frac{\int (\beta - \alpha)^{k} (\beta+\alpha) \,  \tilde{R}(d\beta \vert \alpha) + c\int (\beta - \alpha)^{k} \,  \tilde{R}(d\beta \vert \alpha) }{\int (\beta+\alpha+c) \,  \tilde{R}(d\beta \vert \alpha)} \, .
\]
We know from Appendix A of \citet{Eaton:2008} that
\begin{align}
\int \beta \,  \tilde{R}(d\beta \vert \alpha) &= \alpha + \phi_{1}(\alpha)\, , \label{1st.moment} \\
\int \beta^{2} \,  \tilde{R}(d\beta \vert \alpha) &= \alpha^{2} + 8 \alpha + \phi_{2}(\alpha) \, \label{2nd.moment} \text{, and}  \\
\int \beta^{3} \,  \tilde{R}(d\beta \vert \alpha) &= \alpha^{3} +  24 \alpha^{2} + \phi_{3}(\alpha) \, . \label{3rd.moment}
\end{align}
where, as $\alpha \to \infty$, $\phi_{1}(\alpha) = O(\alpha^{-1})$, $\phi_{2}(\alpha)=O(1)$, and $\phi_{3}(\alpha)=O(\alpha)$.
A similar argument shows that
\begin{equation}
\int \beta^{4} \,  \tilde{R}(d\beta \vert \alpha) = \alpha^{4} + 48 \alpha^{3} + \phi_{4}(\alpha) \label{4th.moment}
\end{equation}
where $\phi_{4}(\alpha) = O(\alpha^{2})$ as $\alpha \to \infty$.  

Let $g_{0}(z) = (a+z)^{-p/2}$, and let $t_{k}(y) = \E[ g_{0}(U)U^{k} \given y ]$ with $U \sim \chi^{2}_{p}(y)$.
Proposition A.2 of \citet{Eaton:2008} establishes that if $Y \sim \chi^{2}_{p}(\alpha)$, then
\[ 
\int \beta^{k} \tilde{R}(d\beta \given \alpha) = \E \left[ \frac{t_{k}(Y)}{t_{0}(Y)} \Big\given \alpha \right]  < \infty \; .
\]
Let 
\[ w_{k}(n) = \int_{0}^{\infty} \frac{g_{0}(z/2)^{n+p/2+k-1}e^{-z/2}}{2\Gamma(n+p/2)} \, dz   \]
and note that $\E[ g_{0}(U)U^{k} \given y ] = 2^{k} \E[ w_{k}(N) \given y ]$ where $N \given y \sim \text{Poisson}(y/2)$.
This last equality follows from expressing $U$ as a Poisson mixture of $\chi^{2}_{p}$ random variables.
From the definition of $t_{k}$ and $w_{k}$, we have that $r_{k}(y) = 2^{k} \E[ w_{k}(N) \given y ]$.
Therefore, 
\[
 \int \beta^{k} \tilde{R}(d\beta \given \alpha) = 2^{k} \E \left[ \frac{\E (w_{k}(N) \given Y)}{\E (w_{0}(N) \given Y)} \Big\given \alpha \right] .
 \]
Proposition A.7 of \citet{Eaton:2008} establishes that
\[ 
\frac{\E [ w_{k}(N) \given y ] }{\E [ w_{k-1}(N) \given y ] } = \frac{y}{2} + 2(k  -1) - \frac{p}{2} + \psi_{k}(y)  
\]
where $\lvert \psi_{k}(y) \rvert \leq d_{k}/y$ for some finite positive constant $d_{k}$.
This allows us to evaluate the right hand side of as 
\[   
\int \beta^{k} \tilde{R}(d\beta \given \alpha) = 2^{k} \E \left[ \frac{\E (w_{k}(N) \given Y)}{\E (w_{k-1}(N) \given Y)} \dotsm \frac{\E (w_{1}(N) \given Y)}{\E (w_{0}(N) \given Y)} \Big\given \alpha \right] . 
\]
Furthermore, $\E[ \psi_{k}(Y) \given \alpha] = O(\alpha^{-1})$ by Proposition A.8 of \citet{Eaton:2008}.
Note that the third moment of a non-central $\chi^{2}_{p}(\alpha)$ is $\alpha^{3} + O(\alpha^{2})$ and the fourth moment is
\[ 
\alpha^{4} + 24\alpha^{3} + 4p\alpha^{3} + O(\alpha^{2}) \; .
\]
Combining these results leads to \eqref{4th.moment}.

We now begin to find expressions for the transitional moments $m_{k}$ in terms of equations \eqref{1st.moment}, \eqref{2nd.moment}, \eqref{3rd.moment}, and \eqref{4th.moment} as $\alpha$ becomes large.
First, note that if $0 < c < \infty$, then
\[   
c \int (\beta - \alpha)^{2} \, \tilde{R}(d\beta \vert \alpha) = O(\alpha) \quad \text{and} \quad c \int (\beta - \alpha)^{3} \, \tilde{R}(d\beta \vert \alpha) = O(\alpha^{2}) .
\]
Express the first transitional moment as
\begin{align*}
m_{1}(\alpha) 
&= \frac{\int (\beta - \alpha)(\beta + \alpha) \, \tilde{R}(d\beta \vert \alpha) + c \int (\beta - \alpha) \, \tilde{R}(d\beta \vert \alpha) }{\int (\beta + \alpha + c) \, \tilde{R}(d\beta \vert \alpha)} \\
&=  \frac{\int \beta^{2} \, \tilde{R}(d\beta \vert \alpha) - \alpha^{2} + c\int \beta \, \tilde{R}(d\beta \vert \alpha) - c\alpha }{\int (\beta + \alpha + c) \, \tilde{R}(d\beta \vert \alpha )}\\
&= \frac{8\alpha + O(1)}{\int (\beta + \alpha +c) \, \tilde{R}(d\beta \vert \alpha )}.
\end{align*}
Since $(\beta - \alpha)^{2} (\beta+\alpha) = \beta^{3} -\beta^{2}\alpha - \beta \alpha^{2} + \alpha^{3}$, we have that
\begin{align*}
\int (\beta - \alpha)^{2} (\beta+\alpha) \,  \tilde{R}(d\beta \vert \alpha) 
&= \alpha^{3} + 24  \alpha^{2} + O(\alpha)  - \alpha^{3} -  8
\alpha^{2} - O(\alpha) - \alpha^{3} - O(\alpha)   + \alpha^{3} \\
&= 16 \alpha^{2} + O(\alpha) \, .
\end{align*}
Express the second transitional moment as
\begin{align*}
m_{2}(\alpha) 
&= \frac{\int (\beta - \alpha)^{2}(\beta + \alpha) \, \tilde{R}(d\beta \vert \alpha) + c \int (\beta - \alpha)^{2} \, \tilde{R}(d\beta \vert \alpha) }{\int (\beta + \alpha + c) \, \tilde{R}(d\beta \vert \alpha)} \\
&= \frac{16\alpha^{2} + O(\alpha)}{\int (\beta + \alpha +c) \, \tilde{R}(d\beta \vert \alpha )}.
\end{align*}
Since $(\beta - \alpha)^{3} (\beta+\alpha) = \beta^{4} - 2\beta^{3}\alpha  + 2\beta\alpha^{3} - \alpha^{4}$, we have that 
\begin{align*}
\int (\beta - \alpha)^{3} (\beta+\alpha) \,  \tilde{R}(d\beta \vert \alpha) 
&= \alpha^{4} + 48 \alpha^{3} + O(\alpha^{2}) - 2\alpha^{4} - 48 \alpha^{3} - O(\alpha^{2}) \\
&\quad + 2\alpha^{4} + O(\alpha) - \alpha^{4} \\
&= O(\alpha^{2}) .
\end{align*}
Since our expressions for the transitional moments all share a common denominator, the ratio $m_{3}(\alpha)/m_{2}(\alpha)$ may be evaluated as
\begin{align*}
\frac{m_{3}(\alpha)}{m_{2}(\alpha )} 
&= \frac{\int (\beta - \alpha)^{3} (\beta+\alpha) \, \tilde{R}(d\beta \vert \alpha) + c\int(\beta-\alpha)^{3} \, \tilde{R}(d\beta \vert \alpha)}{\int (\beta - \alpha)^{2} (\beta+\alpha) \, \tilde{R}(d\beta \vert \alpha) + c\int (\beta - \alpha)^{2}\, \tilde{R}(d\beta \vert \alpha) } \\
&= \frac{O(\alpha^{2})}{16\alpha^{2}+O(\alpha)} \\
&= O(1) .
\end{align*}

\bibliographystyle{plain}

\end{document}